 \documentclass{amsart}
 
 \usepackage{amscd}
 \usepackage[curve,matrix,arrow]{xy}

\usepackage{amsthm}

\newtheorem{lem}[equation]{Lemma}
\newtheorem{cor}[equation]{Corollary}
\newtheorem{prop}[equation]{Proposition}

\newtheorem{MT}[equation]{Main Theorem}

\newtheorem*{thm*}{Theorem}
\newtheorem*{prop*}{Proposition}
\newtheorem*{cor*}{Corollary}
\newtheorem*{lem*}{Lemma}
\newtheorem*{MT*}{Main Theorem}

\theoremstyle{definition} %

\newtheorem*{defn*}{Definition}

\newtheorem{eg}[equation]{Example}

\theoremstyle{remark} %
\newtheorem{rmk}[equation]{Remark}

\newtheorem*{rmk*}{Remark}
\newtheorem*{rmks*}{Remarks}

\newtheoremstyle{exercise}
  {3pt}
  {3pt}
  {\small}
  {\parindent}
  {\sc\small}
  {.}
  {.5em}
   {}     
  {}

\theoremstyle{exercise}

\makeatletter
{\renewcommand{\theequation}{#1}}%
{\renewcommand{\theequation}{\arabic{equation}}\addtocounter{equation}{-1}\global\@ignoretrue}
\makeatother

\makeatletter
{\renewcommand{\theequation}{#1}\begin{eqnarray}}%
{\end{eqnarray}\renewcommand{\theequation}{\arabic{equation}}\addtocounter{equation}{-1}\global\@ignoretrue}
\makeatother

\makeatletter
\newenvironment{borel}[1]%
{\smallskip \refstepcounter{equation}\noindent{\textbf \theequation. }{{\textbf{#1.}}}}%
{\smallskip \global\@ignoretrue}
\makeatother

\makeatletter
{\smallskip \refstepcounter{equation}{\sc \theequation}{\sc (#1).}}%
{\smallskip \global\@ignoretrue}
\makeatother

\makeatletter
{\smallskip \refstepcounter{equation}\noindent{\sc \theequation.}{\sl{ #1.}}}%
{\smallskip \global\@ignoretrue}
\makeatother

\makeatletter
\newenvironment{borel*}%
{\smallskip \refstepcounter{equation}\noindent{\textbf \theequation.}}%
{\global\@ignoretrue}
\makeatother

\newcommand{\flist}[1]{\hangindent\leftmargini\textup{(1)}\hskip\labelsep {#1}%
\begin{enumerate}%
\setcounter{enumi}{1}%
}

\newcommand{\ot}{\otimes}

\newcommand{\Q}{{\mathbb{Q}}}        
\newcommand{\R}{{\mathbb{R}}}        
\newcommand{\Z}{{\mathbb{Z}}}        

\newcommand{\QZt}{\mathbb{Q}/\mathbb{Z}(2)}

\newcommand{\Zm}[1]{\Z/{#1}\Z}

\newcommand{\La}{\Lambda}
\newcommand{\la}{\lambda}
\newcommand{\D}{\Delta}

\newcommand{\G}{{\Gamma}}       

\newcommand{\oddots}{{\mathinner{\mkern1mu\raise1pt\vbox{\kern7pt\hbox{.}}\mkern2mu\raise4pt\hbox{.}\mkern2mu\raise7pt\hbox{.}\mkern1mu}}}

\newcommand{\s}{\sigma}

\newcommand{\ksep}{k_{{\mathrm{sep}}}}

\newcommand{\supast}[1]{{{#1}^\times}}
\newcommand{\xsq}[1]{{{#1}^{\times2}}}
\newcommand{\sq}[1]{{\supast{#1}/\xsq{#1}}}

\newcommand{\kx}{k^\times}
\newcommand{\ksq}{\sq{k}}

\newcommand{\qform}[1]{{\langle{#1}\rangle}}                   

\newcommand{\basemu}{\boldsymbol{\mu}}
\newcommand{\mmu}[1]{\basemu_{#1}}     
\newcommand{\mmut}[1]{{\basemu^{\otimes 2}_{#1}}}     

\DeclareMathOperator{\Spin}{Spin}           

\newcommand{\Ga}{\mathbb{G}_a}
\newcommand{\Gm}{\mathbb{G}_m}

\newcommand{\iiiD}{^3\!D_4}
\newcommand{\viD}{^6\!D_4}
\newcommand{\oD}{^1\!D_4}

\newcommand{\iiD}{^2\!D_4}

\DeclareMathOperator{\Gal}{Gal}

\DeclareMathOperator{\Int}{Int}
\DeclareMathOperator{\chr}{char}

\DeclareMathOperator{\aut}{Aut}

\newcommand{\Hom}{{\mathrm{Hom}}}

\newcommand{\ra}{\rightarrow}

 \numberwithin{equation}{section}
 
 \newcommand{\mud}[1]{\boldsymbol{\mu}_{{#1}}^{\otimes 2}}
 
 \newcommand{\darkrad}{0.115}
\newcommand{\lrad}{0.25}

\newcommand{\even}{{{\mathrm{even}}}}
\newcommand{\odd}{{{\mathrm{odd}}}}

\newcommand{\cprod}{{\mathinner{\mkern2mu\raise2.5pt\hbox{.}}}}
\newcommand{\cp}{\cprod}

\newcommand{\sep}{{\mathrm{sep}}}

\newcommand{\ks}{k_\sep}

\newcommand{\Tbar}{\overline{T}}

\newcommand{\Mbar}{\overline{M}}

\DeclareMathOperator{\Inv}{Inv}

\newcommand{\kxsq}{k^{\times 2}}

\newcommand{\co}{{\vee}}
\newcommand{\coroot}{\La_r^\co}

\newcommand{\coweight}{\La^\co}

\newcommand{\chom}{{\check{\omega}}}
\newcommand{\och}{\chom}

\newcommand{\ach}{{\check{\alpha}}}
\newcommand{\Pch}{\Phi^\co}

\newcommand{\Gbar}{\overline{G}}

\newcommand{\cF}{{\mathcal F}}
 
\begin{document}

\title{Restricting the Rost Invariant to the center}
\author{Skip Garibaldi}
\address{Department of Mathematics \& Computer Science, Emory University, Atlanta, GA 30322, USA}
\email{skip@member.ams.org}
\urladdr{http://www.mathcs.emory.edu/{\textasciitilde}skip/}

\author{Anne Qu\'eguiner-Mathieu} 
\address{Laboratoire Analyse, G\'eom\'etrie \& Applications\\
UMR CNRS 7539 -- Institut Galil\'ee\\
Universit\'e Paris 13\\
93430 Villetaneuse\\
France}
\email{queguin@math.univ-paris13.fr}
\urladdr{http://www-math.univ-paris13.fr/{\textasciitilde}queguin/}

\begin{abstract} For simple simply connected algebraic groups of
  classical type, Merkurjev, Parimala, and Tignol gave a formula for
  the restriction of the Rost invariant to torsors induced from the
  center of the group. 
We complete their results by proving formulas for exceptional groups. 
Our method is somewhat different and recovers also their formula for
  classical groups. 
\end{abstract}

\date{\today}
\maketitle

\section{Introduction}
In the 1990's, M. Rost proved that the group of degree $3$ normalized invariants of
an absolutely simple simply connected algebraic group $G$ over $k$
--- that is, the group of natural transformations of the Galois cohomology functors  
\[
H^1(\star,G)\ra H^3(\star,\QZt)
\]
--- is a cyclic group with a canonical generator.  This canonical generator is known as the \emph{Rost invariant}.  Roughly speaking, it is the ``first" nonzero invariant, in that 
there are no non-zero normalized invariants $H^1(\star, G) \ra H^d(\star, \Q/\Z(d - 1))$ for $d < 3$ \cite[\S31]{KMRT}.  

In general, there is no explicit formula for computing the Rost invariant.  However, in \cite{MPT}, Merkurjev, Parimala and Tignol gave a nice
description --- for $G$ of classical type --- of the restriction of
the Rost invariant to $H^1(\star,Z)$, where $Z$ is the center of $G$. The purpose of this paper is to complete their results by providing an analogous
description for exceptional groups (see Th.~\ref{mainthm.thm}
for a precise statement for all types of groups). 

Part of our proof is borrowed from \cite{MPT}. Precisely, their
Corollaries 1.2 to 1.6, which are stated for a general cycle module, 
show in our context that every group invariant 
\[
H^1(\star,Z)\ra H^3(\star,\QZt)
\]
--- i.e., every invariant such that the map $H^1(K, Z) \ra H^3(K, \QZt)$ is a group homomorphism for every extension $K/k$ --- 
is given by some cup-product with a class $t\in
H^2(k,Z)$, see Prop.~\ref{invZ.prop} and \S\ref{invZ.section} below. 
(Note that this class $t$ does depend on a choice regarding
the cup product, see Remark~\ref{Deven.rmk}(iii).)  

It remains to compute the class 
$t_{R,G}\in H^2(k,Z)$ associated with composition 
\begin{equation} \label{hom}
\begin{CD}
H^1(\star, Z) @>>> H^1(\star, G) @>{r_G}>> H^3(\star, \QZt),
\end{CD}
\end{equation}
where the first map is induced by the inclusion of $Z$ in $G$ and 
$r_G$ is the Rost invariant of $G$. 
Note that, even though it is not obvious from its definition, 
this composition actually is a group invariant, see \cite[7.1]{G:rinv} or \cite[Cor.~1.8]{MPT}.  The Rost
invariant is canonically determined; nevertheless, practically
speaking, one only knows the group it generates.  Consequently, we will determine the subgroup $\qform{t_{R,G}}$ of $H^2(k, Z)$ generated by $t_{R,G}$.

For any semisimple group $G$, one has the
Tits class $t_G$, an element of the group $H^2(k,Z)$.
It follows from~\cite{MPT} and from the present paper that
\emph{either $t_{R,G}$ is zero or $t_{R,G}$ and $t_G$ generate the same subgroup of $H^2(k, Z)$}, for a well-chosen cup product.  
Note that our argument for determining $\qform{t_{R,G}}$ is quite
different from~\cite{MPT}, in which they use  
concrete interpretations of the classical groups. 
Here, we first reduce to groups whose Tits
index satisfies a certain condition (see \eqref{cond}), 
using an injectivity result which follows from~\cite[Th.~B]{MT}. 
Second, we prove the result for those particular groups, 
by reducing to groups of inner type $A$ for which 
the Rost invariant has a concrete description.
This method applies to all absolutely simple groups except those of outer type $A$, and therefore can be used to recover the results of \cite{MPT} on classical groups (except for their Th.~1.10).  We give the details for types  $B$, $C$, and
$D_{{\mathrm{even}}}$, as well as for the exceptional groups.

\section{Statement of results} \label{results.section} 

Throughout the paper, $G$ denotes an absolutely simple simply
connected algebraic group over a field $k$.  
We assume (except in \S\ref{center.1}) that \emph{the characteristic of $k$ does not divide the exponent $n$ of the center of $G$}.  That is, the characteristic is not 2 for $G$ of type $B$, $C$, $D$, $E_7$; the characteristic is not 3 for $G$ of type $E_6$; and the characteristic does not divide $\ell+1$ for $G$ of type $A_\ell$.  This guarantees that the (scheme-theoretic) center $Z$ of $G$ is smooth.  (See Remark \ref{reduction} for comments on the characteristic hypothesis.) 

The group $H^d(k, \Q/\Z(d - 1))$ is as defined in \cite[App.~A]{MG}; otherwise the notation $H^d(k, A)$ stands for the Galois cohomology group $H^d(\G_k, A(\ksep))$, where $\G_k$ denotes the Galois group of a separable closure $\ksep$ of $k$ over $k$ and $A$ is a smooth algebraic group.  Our hypothesis on the characteristic of $k$ implies that the $n$-torsion of $H^d(k, \Q/\Z(d - 1))$ is naturally identified with $H^d(k, {\boldsymbol{\mu}}^{\otimes{(d - 1)}}_n)$ for $d = 2, 3$.
When $G$ is not of type $A_\ell$ with $\ell \ge 3$ --- e.g., in all
cases studied in detail below --- $n$ is 2 or 3, hence $\mud{n}$ is
isomorphic to $\Zm{n}$ (see \cite[Ex.~11, p.~444]{KMRT}). 

For any functor $\cF$ from the category of field extensions of $k$ to
the category of groups, we let $\Inv^3(\cF)$ be the collection of
group invariants of $\cF$ with values in $H^3(\star,\QZt)$; it is an abelian group.
If $T$ is a quasi-trivial torus, i.e., $T=R_{E/k}(\Gm)$ for some 
\'etale algebra $E/k$, it follows from Th.~1.1 of \cite{MPT}
that $H^2(k, T)$ is isomorphic to 
$\Inv^3(T)$. 
Precisely, the invariant $\alpha^E(X)$ associated to a cohomology
class $X\in H^2(k, T)$ is defined by 
$\alpha^E(X)(y)=N_{E_K/K}(y \cp X_K)$ for every field extension $K/k$ 
and $y \in T(K)=E_K^\times$. The cup product
appearing in this formula is given by the module structure on
$H^\star(\Q/\Z(-1))$ over the Milnor $K$-ring 
(see~\cite[Appendix A]{MG} for a definition). 
Using this, we may prove the following: 

\begin{prop} \label{invZ.prop}
The groups $\Inv^3(H^1(\star,Z))$ and $H^2(k,Z)$ 
are (non-canonically) isomorphic. 
\end{prop}

While proving this proposition in section~\ref{invZ.section}, we will
exhibit particular such isomorphisms, which admit the
following explicit description: for any cohomology class $X\in
H^2(k,Z)$, the corresponding invariant is given by 
\[
y\in H^1(K,Z)\mapsto y \cp X_K\in H^3(K,\QZt).
\] 
Since it is a group invariant, its image is contained in the $n$-torsion, $H^3(K, \mud{n})$, of $H^3(K, \QZt)$.  Further, 
the cup product is
induced by a bilinear form $Z(\ks)\times Z(\ks)\ra \mud{n}$, which
has to be specified.

We denote by $\mmu{n}$ the
algebraic group of $n$th roots of unity and by ${\mmu{n}}_{[E]}$ 
the kernel of the norm map 
$N_{E/k}:\,R_{E/k}(\mmu{n})\ra \mmu{n}$, for any quadratic \'etale algebra
$E/k$.  
As recalled in~\cite{MPT}, one deduces from the classification of 
absolutely simple simply connected groups that if $G$ is classical, 
then $Z$ is one of the following groups:
$\mmu{n}$, $R_{E/k}(\mmu{2})$, and ${\mmu{n}}_{[E]}$, where $E$ is quadratic \'etale
over $k$. 
If we wish to consider also exceptional groups, we
need to add the centers of trialitarian $D_4$ groups, which are
isomorphic to the kernel, now denoted by 
$R^1_{E/k}(\mmu{2})$, of the norm map 
$N_{E/k}:\,R_{E/k}(\mmu{2})\ra \mmu{2},$ 
where $E$ is cubic \'etale
over $k$. 

In most cases, namely when $Z$ is $\mmu{n}$ or ${\mmu{n}}_{[E]}$, there is a natural
  bilinear map 
\[
Z(\ks)\times Z(\ks)\ra\mud{n}.
\] 
The only groups for which the cup-product has to be defined carefully
are groups of type $D_\ell$ with $\ell$ even.
Their center is---over $\ksep$---isomorphic to $\mmu2 \times \mmu2$.  Rather than fix an identification between these two groups, we note that the fundamental weights of the root system give characters $\omega_1, \omega_{\ell-1}, \omega_\ell$ whose restriction to $Z$ are the three nonzero homomorphisms $Z \ra \mmu2$.  We consider the cup
product induced by the bilinear map with values in $\mud{2} = \Zm2$:
\begin{equation} \label{Deven.cup}
(x, y) \mapsto 
\begin{cases}
\omega_{\ell-1}(x) \ot \omega_\ell(y) + \omega_{\ell}(x) \ot \omega_{\ell-1}(y) & \text{if $\ell \equiv 0 \pmod{4}$,}\\
\omega_{\ell-1}(x) \ot \omega_{\ell-1}(y) + \omega_{\ell}(x) \ot \omega_{\ell}(y) & \text{if $\ell \equiv 2 \pmod{4}$.}
\end{cases}
\end{equation}

Since $\omega_1 + \omega_{\ell-1}+\omega_\ell = 0$ as characters, when
$\ell\equiv 0 \pmod{4}$, this cup product can be rewritten as
\[
(x, y) \mapsto \omega_1(x) \ot \omega_1(y) + \omega_{\ell-1}(x) \ot \omega_{\ell-1}(y) + \omega_\ell(x) \ot \omega_\ell(y).
\]

With this in hand, we may summarize the results of~\cite{MPT} and the
present paper in the following theorem:\footnote{There is a typo in Th.~1.13 of \cite{MPT}, which addresses
  the $C_\ell$ case: the words ``odd'' 
and ``even" should be interchanged.  
For keeping the $C_\odd$ and $C_\even$ cases straight, we find it
helpful to 
recall that $B_2 = C_2$, cf.~\S\ref{C.section} below.}
\begin{MT} \label{mainthm.thm}
If $G$ is of type $A$, $C_\ell$ ($\ell$ odd), $D$, $E_6$ or $E_7$, then the
composition \eqref{hom} and 
the cup product with the Tits class of $G$ generate the same subgroup
of 
$\Inv^3(H^1(\star,Z))$.  Otherwise the composition \eqref{hom} is zero.
\end{MT}

\begin{eg}\label{typeA.rmk} 
For groups of inner type $A$ --- that is, when $G$ is $SL(A)$ for some central
simple algebra $A$ over $k$ --- the theorem asserts that \eqref{hom} and the cup product with the Brauer class $[A]$ of $A$ generate the same subgroup of $\Inv^3(H^1(\star, Z))$.
This follows from the explicit description of the Rost invariant given
in \cite[p.~138]{MG}. 
This example is a starting point for proving the theorem. 
\end{eg} 

The new cases in the theorem are the groups of type $E_6$, $E_7$, and
trialitarian $D_4$, 
treated below in sections \ref{E6.sec}, \ref{E7.sec}, and \ref{D.section}
respectively.  
The exceptional groups of type $E_8$, $F_4$, and $G_2$ all have
trivial center, so \eqref{hom} 
is automatically zero for those groups.

\begin{rmk} \label{Deven.rmk}
(i) \emph{If the Tits class is zero, then a twisting argument shows directly that the composition \eqref{hom} is zero, regardless of the type of $G$.}  Indeed, since the Tits class is zero, there is a simply connected quasi-split group $G^q$ such that we may identify $G$ with $G^q$ twisted by a 1-cocycle $\alpha \in Z^1(k, G^q)$.  We write $Z^q$ for the center of $G^q$; the identification of $G$ with $G^q_\alpha$ identifies $Z$ with $Z^q$; we write $i$ and $i^q$ for the inclusions of $Z$ and $Z^q$ in $G$ and $G^q$.  The diagram
\[
\begin{CD}
H^1(k, Z) @>{i_\star}>> H^1(k, G) @>{r_G}>> H^3(k, \QZt) \\
@. @V{\tau_\alpha}V{\cong}V @VV{? + r_{G^q}(\alpha)}V \\
H^1(k, Z^q) @>{i^q_\star}>> H^1(k, G^q) @>{r_{G^q}}>> H^3(k, \QZt)
\end{CD}
\]
where $\tau_\alpha$ is the twisting isomorphism, commutes by \cite[p.~76, Lemma 7]{Gille:inv} or
\cite[Cor.~1.8]{MPT}.  Fix $\zeta \in Z^1(k, Z)$.  The image of
$i_\star \zeta$ under $\tau_\alpha$ is the class in $H^1(k, G^q)$ of the 1-cocycle $\s \mapsto \zeta_\s \alpha_\s$ for $\s \in \G_k$.  But this is just $\zeta \cdot \alpha \in H^1(k, G^q)$, where the $\cdot$ represents the usual action of $H^1(k, Z^q)$ on $H^1(k, G^q)$.  So the image of $i_\star \zeta$ going counterclockwise around the diagram is 
\[
r_{G^q}(\zeta \cdot \alpha) = r_{G^q}(i^q_\star \zeta) + r_{G^q}(\alpha)
\]
by \cite[7.1]{G:rinv}.  But $G^q$ is quasi-split, so it has a quasi-trivial maximal torus $T^q$ that necessarily contains $Z^q$, hence $i^q_\star$ factors through the zero group $H^1(k, T^q)$.    We conclude that the image of $i_\star \zeta$ in the lower right corner is $r_{G^q}(\alpha)$.

On the other hand, going clockwise around the diagram, we find that the image of $i_\star \zeta$ in the lower right corner is $r_G(i_\star \zeta) + r_{G^q}(\alpha)$.  The commutativity of the diagram proves the claim.

\medskip
(ii) As opposed to (i), note:
For groups of type $B_\ell$, i.e., for $G=\Spin(q)$ for
some $(2\ell+1)$-dimensional quadratic form $q$ over $k$, the composition \eqref{hom} is always zero, while the Tits class of $G$,
which is the Brauer class of the even Clifford algebra of $q$, can
be non-zero. 

\medskip
(iii) In the $D_{\ell}$ case with $\ell$ even, the bilinear form 
we use to compute the cup product does depend 
on whether $\ell/2$ is even or odd.
If we keep the same cup product in both cases, then composition
\eqref{hom} 
coincides alternately with the cup product
with $t_G$ and its conjugate. 

\medskip
(iv) The bilinear maps in \eqref{Deven.cup} can be viewed as
 symmetric bilinear forms on a vector space of dimension 2 over the field with
  two elements.  For $\ell \equiv 0 \pmod 4$, it is the wedge product, equivalently, a 
  hyperbolic form.  For $\ell \equiv 2 \pmod 4$, it is the unique (up to isomorphism) metabolic form that is not hyperbolic.
\end{rmk} 

\section{Invariants of $H^1(\star,Z)$} \label{invZ.section}

The purpose of this section is to prove Proposition~\ref{invZ.prop}. 
In order to use the description of invariants of quasi-trivial tori
recalled above, we need to embed $Z$ in such a torus.
More precisely, we are going to consider exact sequences 
\[
\begin{CD}
1 @>>> Z @>{j}>> T @>{g}>> T @>>> 1,
\end{CD}
\]
where $T$ is a quasi-trivial torus, i.e., $T=R_{E/k}(\Gm)$ for some \'etale
$k$-algebra $E$. 
We call such an exact sequence \emph{admissible} if the map $g$ satisfies
the following equation: 
\begin{equation}\label{admissible.cond}
N_{E_K/K}(g_K(y) \cp X_K)=N_{E_K/K}(y \cp H^2(g)(X)_K),
\end{equation}
for every $X\in H^2(k,T)$, field extension $K/k$, and $y\in
T(K)=E_K^\times$. 

Given an admisible exact sequence, we claim it induces an isomorphism 
between $\Inv^3(H^1(\star,Z))$ and $H^2(k,Z)$, as in
Prop.~\ref{invZ.prop}.
Indeed we can apply \cite[Lemma 3.1]{MPT}: Since the sequence
\[
\begin{CD}
T(K) @>{g_K}>> T(K) @>{\varphi}>> H^1(K,Z) @>>> 1. 
\end{CD}
\]
is exact for every extension $K/k$, the top row of the diagram
\begin{equation} \label{iso.diag}
\xymatrix{
1 \ar[r] & \Inv^3(H^1(\star,Z)) \ar[r]^{\varphi^\star} & \Inv^3(T) \ar[r]^{g^\star} & \Inv^3(T)\\
1 \ar[r] & H^2(k,Z) \ar@{-->}[u]^{\beta} \ar[r]^{H^2(j)} & H^2(k,T) \ar[u]^{\alpha^E}_{\cong} \ar[r]^{H^2(g)} & H^2(k,T) \ar[u]^{\alpha^E}_{\cong}
}
\end{equation}
is also exact, where $\alpha^E$ is the isomorphism described in \S\ref{results.section}.  Since $T$ is quasi-trivial, the group $H^1(k, T)$ is zero and the bottom row of the diagram is also exact.  An easy computation shows that under condition~\eqref{admissible.cond}
the right-hand box is commutative. Since $\alpha^E$ is an isomorphism, this enable us to 
identify the kernels of both lines, producing an isomorphism $\beta$ as in the diagram.

To prove Prop.~\ref{invZ.prop}, it only remains to show that
there exists an admissible exact sequence for the center $Z$ of a given $G$. 
This is done in~\cite{MPT} in sections 3.1 to 3.4 
when $Z$ is isomorphic to $\mmu{n}$, $R_{E/k}(\mmu{n})$, or ${\mmu{n}}_{[E]}$ 
(see also sequences (3), (6), and (8)). The only
remaining case is the center of a trialitarian
$D_4$. For such a group, there exists a cubic field extension 
$E/k$ and an embedding $j:\,Z\ra R_{E/k}(\Gm)$ which identifies $Z$ with
the kernel $R^1_{E/k}(\Gm)$ of the norm map. 
Over
$\ksep$, it is given by  
$j_\sep(x)=(\omega_1(x),\omega_{\ell-1}(x),\omega_\ell(x))\in (\ksep^\times)^3$.

Let us define $g\!:R_{E/k}(\Gm)\ra R_{E/k}(\Gm)$ 
by 
$g_K(x)=N_{E_K/K}(x)\,x^{-2}$, for any field $K/k$
and any $x\in R_{E/k}(\Gm)(K)={E}_K^\times$, and consider the sequence 
\begin{equation} \label{D4seq.equation} 
\begin{CD} 
1 @>>> Z @>{j}>> R_{E/k}(\Gm) @>{g}>> R_{E/k}(\Gm) @>>> 1
\end{CD} 
\end{equation} 
\begin{lem} 
Sequence \ref{D4seq.equation} is exact and admissible.  
\end{lem}

\begin{proof}
An element $x$ in the kernel of $g_K$ satisfies 
$N_{E_K/K}(x)=x^{2}$. 
Since $E$ is cubic, taking the norm on both sides gives 
$1 = N_{E_K/K}(x)=x^2$, and this
proves $x\in Z(K)$. 
The reader may easily check surjectivity;  admissibility
condition~\ref{admissible.cond} follows from the projection formula, 
as in~\cite[3.3]{MPT}. 
\end{proof}
Proposition~\ref{invZ.prop} is now proved. 

Note that, if you allow $E/k$ to be a cubic \'etale algebra rather
than a field extension, then the sequence~(\ref{D4seq.equation}) actually is an
admissible exact sequence for any group of type $D_\ell$ with
$\ell\equiv 0\mod 4$. 
As we will see in the next paragraph, the corresponding isomorphism 
between $\Inv^3(H^1(\star,Z))$ and $H^2(k,Z)$ is not the same as the
one induced by sequence (3) of \cite{MPT} (see also
Remark~\ref{Deven.rmk}(iii)). 

Hence we have to specify 
which exact sequence we use: 
from now on, we will always consider the isomorphism induced by the
following admissible exact sequence, depending on $Z$ and $G$:
\begin{itemize}
\item If $Z = \mmu{n}$, we take
\[
\begin{CD} 
1 @>>> \mmu{n} @>>> \Gm @>{n}>> \Gm @>>> 1.
\end{CD}
\]

\item If $Z = {\mmu{n}}_{[E]}$ with $n$ odd (respectively
    even), we take sequence (6) (resp., (8)) 
of~\cite{MPT} 

\item If $G$ is of type $D_{\ell}$ with $\ell\equiv 2\mod 4$, we take sequence (3) of \cite{MPT}.

\item If $G$ is of type $D_{\ell}$ with $\ell\equiv 0\mod 4$, we take sequence \eqref{D4seq.equation}.
\end{itemize}
In each of these cases, the isomorphism $\beta$ from \eqref{iso.diag}
is given by 
\begin{equation} \label{cup.eqn}
\beta(t)(x)=x\cp t_K,
\end{equation}
where the cup product is induced by
the bilinear map specified in \S~\ref{results.section}.   (This is the description of the isomorphism $\beta$ given just after the statement of Prop.~\ref{invZ.prop}.)
This fact can be proved in each case by some explicit computation,
going through the diagram. 
Let us do it for groups of type $D_{\ell}$ with $\ell\equiv 0\mod 4$,
and sequence \eqref{D4seq.equation}.

For any $t\in H^2(k,Z)$, the invariant $\beta(t)$ is characterized by 
\begin{equation} \label{cup.1}
\beta(t)(\varphi(x))=\alpha^E(H^2(j)(t))(x)
=N_{E_K/K}(x\cp H^2(j)(t_K))
\end{equation}
for any $x\in T(K)$. 
We want to prove this is equal to $\varphi(x)\cp t_K$, where the cup
product is induced by the bilinear map 
\[
(x, y) \mapsto \omega_1(x) \ot \omega_1(y) + \omega_{\ell-1}(x) \ot \omega_{\ell-1}(y) + \omega_\ell(x) \ot \omega_\ell(y).
\]
Since $Z$ has exponent $2$, the map $j$ factors as 
\[\begin{CD} 
Z@>{j_0}>>R_{E/k}(\mmu{2})@>>>R_{E/k}(\Gm),\end{CD}\]
and ${j_0}_\sep$, as $j_\sep$, coincides with the triple
$(\omega_1,\omega_{\ell-1},\omega_\ell)$. 
Hence, we have 
\begin{equation} \label{cup.2}
\varphi(x)\cp t_K=N_{E_K/K}(H^1(j_0)(\varphi(x))\cp H^2(j_0)(t_K)).
\end{equation}
To compare \eqref{cup.1} and \eqref{cup.2}, we use the following commutative
diagram
\[\xymatrix{
1 \ar[r] & Z \ar[d]^{j_0}\ar[r]^{j} &R_{E/k}(\Gm)\ar@{=}[d] \ar[r]^{g} & R_{E/k}(\Gm)\ar[d]^{h}\ar[r]&1\\
1 \ar[r] & R_{E/k}(\mmu{2}) \ar[r] & R_{E/k}(\Gm)  \ar[r]^{2} & R_{E/k}(\Gm) \ar[r]&1
}\]
where the map $h$ is defined by $h_K(x)=N_{E_K/K}(x)x^{-1}=xg_K(x)$ for any
$x\in R_{E/k}(\Gm)(K)=E_K^\times$. 
It induces a commutative square 
\[\xymatrix{
E_K^\times\ar[r]^{\varphi}\ar[d]^{h_K}&H^1(K,Z)\ar[d]^{H^1(j_0)}\\
E_K^\times\ar[r]&H^1(K,R_{E/k}(\mmu{2}))\\
}\]
Hence $H^1(j_0)(\varphi(x))=(h_K(x))_2=(x)_2+(g_K(x))_2$, 
and
\[
\begin{array}{r@{=}l}
\multicolumn{2}{l}{N_{E_K/K}(H^1(j_0)(\varphi(x))\cp H^2(j_0)(t_K))} \\
\parbox{1in}{\hfill}&\, N_{E_K/K}((x)_2\cp
H^2(j_0)(t_K))+N_{E_K/K}((g_K(x))_2\cp H^2(j_0)(t_K))\\
&\, N_{E_K/K}(x\cp H^2(j)(t_K))+N_{E_K/K}(g_K(x)\cp H^2(j)(t_K)),
\end{array}
\]
where the cup product in the last line is again given by the
module structure on $H^\star(\Q/\Z(-1))$ over the Milnor $K$-ring (see
\cite[Appendix A]{MG}). 
By the admissibility condition~(\ref{admissible.cond}), the second term
in the sum is 
$N_{E_K/K}(x\cp H^2(g)\circ
H^2(j)(t_K))=0$, so the expressions in \eqref{cup.1} and \eqref{cup.2} are equal.  This finishes the proof of \eqref{cup.eqn}.

\section{Reduction to groups having a particular Tits index}  \label{Titsindex.section} 

For each group $G$, we let $t_{R,G}$ be the 
class in $H^2(k,Z)$ corresponding to the
composition~\eqref{hom} under the isomorphism specified in the previous
section. 
The main theorem asserts that this class is zero or generates the same subgroup of $H^2(k, Z)$ as the Tits class $t_G$, depending on the type of the group $G$ we
started with. 
A case by case proof will be given in
sections~\ref{B.section} to \ref{D.section}. 
First, we prove some general facts, on which our strategy is based.

Given $G$, we 
fix a maximal $k$-split torus $S$, a maximal $k$-torus $T$ containing
 $S$, and a set of simple roots 
$\D$ of $G$ with respect to $T$. 
Recall the \emph{Tits index of $G$} as defined in \cite{Ti:Cl}.  It is
the data of the Dynkin diagram of $G$ together with the action of the Galois
group $\Gal(\ksep/k)$ on $\D$ 
via the $*$-action, and the set $\D_0$ of those $\alpha \in \D$ that
vanish on $S$.  
Notations for $\D_0$ vary, but following \cite{Ti:Cl}, we circle 
a vertex if it belongs to 
$\D \setminus \D_0$, and circle together vertices that are in the
same Galois orbit. 

Let $\D_r$ be the subset of $\D$
consisting of those simple roots 
such that the corresponding fundamental weight belongs to the root
lattice.  
We consider the following condition on $G$:
\begin{equation} \label{cond}
\parbox{4in}{\emph{No $\alpha \in \D_r$ vanishes on $S$.}}
\end{equation}
In terms of Tits indices, this amounts to
\begin{equation} \label{cond.ind}
\parbox{4in}{\emph{Every vertex in $\D_r$ is circled in the Tits index of $G$.}}
\end{equation}
Dynkin diagrams for $E_7$, $E_6$, and $D_\even$ with the vertices in 
$\D_r$ circled can be found in sections \ref{E7.sec}, \ref{E6.sec}, and 
\ref{D.section} below.  

\begin{borel*} \label{Titsindex.prop}
We now give a method that reduces the proof of the main theorem \ref{mainthm.thm} to the case where $G$ satisfies condition \eqref{cond.ind}.  Let $G$ be a simply connected absolutely almost simple group over $k$ and further suppose that $G$ has inner type, i.e., that the absolute Galois group acts trivially on the Dynkin diagram of $G$.  
Because $G$ has inner type, the center $Z$ of $G$ is isomorphic to
$\mmu{n}$ or $\mmu{2} \times \mmu{2}$, and in \S\ref{invZ.section} we
fixed an injection $j \!: Z \ra \Gm^{\times s}$ for $s\leq 3$.  

Consulting the tables in \cite{Bou:g4}, we observe that every element of $\D_r$ is 
fixed by every automorphism of the Dynkin diagram, hence is fixed by the $*$-action.  It follows that 
the variety of parabolic subgroups of $G_{\sep}$ of type $\D \setminus \D_r$
(in the notation of \cite[4.2, 5.12]{BoTi} or \cite[p.~33]{MT})
is defined over $k$, see \cite[5.24]{BoTi} or \cite[8.4]{BoSp2}.
This variety has a point over its function field $F$,
hence $G$ satisfies condition \eqref{cond.ind} \emph{over $F$}.  We have a commutative diagram with exact rows
\[
\begin{CD}
1 @>>> H^2(k, Z) @>{H^2(j)}>> \oplus^s H^2(k, \Gm) \\
@. @VVV @VVV \\
1 @>>> H^2(F, Z) @>{H^2(j)}>> \oplus^s H^2(F, \Gm).
\end{CD}
\]
The kernel of the restriction map $H^2(k, \Gm) \ra H^2(F, \Gm)$
has been computed by Merkurjev and Tignol in~\cite[Th.~B]{MT}. 
Let $\Lambda/\Lambda_r$ be the quotient 
of the weight lattice of $G$ by its root lattice. 
The kernel of the restriction map is generated by the Tits algebras associated with the classes in
$\Lambda/\Lambda_r$ of the fundamental
weights corresponding to those simple roots which belong to $\D_r$. 
But we have precisely chosen $\D_r$ so that these classes are
zero. 
Hence, the restriction map $H^2(k, \Gm) \ra H^2(F, \Gm)$ is injective.  By commutativity of the diagram, the map $H^2(k, Z) \ra H^2(F, Z)$ is injective.

Let $t_{R,G}$ and the Tits class $t_G \in H^2(k, Z)$ be as defined in the introduction.  The isomorphism $\beta$ between $H^2(k, Z)$ and $\Inv^3(H^1(\star, Z))$ from \S\ref{invZ.section} is functorial in $k$, i.e., the image of $t_{R,G}$ in $H^2(F, Z)$ is the element corresponding to the composition \eqref{hom} over $F$.  Therefore, if the composition \eqref{hom} is zero over $F$, then $t_{R,G}$ is killed by $F$, and $t_{R,G}$ is zero in $H^2(k, Z)$, i.e., the composition \eqref{hom} is zero over $k$.  Alternatively, if the restrictions of $t_{R,G}$ and $t_G$ generate the same subgroup of $H^2(F, Z)$, then $t_{R,G}$ and $t_G$ generate the same subgroup of $H^2(k, Z)$.

In this way, we can reduce the proof of the main theorem for groups of inner type to the case of groups of inner type satisfying \eqref{cond.ind}.  This reduction can be generalized to treat also groups of outer type, but we will not use that, so we omit it.  
\end{borel*}

\begin{rmk} \label{reduction}
(i) For groups of type $A$, $\D_r$ is the empty set, so \eqref{cond.ind} vacuously 
holds, the function field $F$ is $k$ itself, and this reduction is useless.

(ii) Our global hypothesis excluding fields of certain characteristics arises from this argument, i.e., from \ref{Titsindex.prop}.  The hypothesis on the characteristic is needed for the identification of $H^2(k, Z)$ with $\Inv^3(H^1(\star, Z))$ from Prop.~\ref{invZ.prop}, because the proof of that proposition used \cite[Th.~1.1]{MPT} to provide the isomorphism $\alpha^E$.  In turn, the proof of the result from \cite{MPT} needs that the $n$-torsion in $H^3(k, \QZt)$ is a cycle module, in particular, one should be able to take residues relative to discrete valuations on $k$.   This is only known for $n$ not divisible by $\chr k$.  

With small changes, our proof for groups of inner type satisfying \eqref{cond.ind} holds without any restriction on the characteristic of $k$.  
 \end{rmk} 

\section{A semisimple subgroup} \label{subgroup.section} 

In this section, we assume that the group $G$ has Tits index
satisfying condition~\eqref{cond.ind}.  We produce a semisimple simply connected subgroup $G'$ of $G$ that contains the center of $G$ and describe how to compute the Rost invariant of some elements of $G$ using $G'$.


\begin{borel}{Description of $G$ by generators and relations} \label{genandrel}
Recall from \cite{St} that --- over $\ksep$ --- $G$ is generated by the images of homomorphisms 
\[
x_\alpha \!:  \Ga \ra G
\]
as $\alpha$ varies over the set of roots $\Phi$ of $G$ with respect to $T$. 
Write 
$\Lambda_r$ and $\Lambda$ for the root and weight lattices of $G$ with respect to $T$.
Since $G$ is simply connected, $\Lambda$ is identified with the
character group $T^\star=\Hom_{\ksep}(T_{\sep},\Gm)$ and $\Lambda_r$ with the character group
$\Tbar^\star$ of the image $\Tbar$ of $T$ in $G/Z$.

Write $T_\star$ for the $\Gamma_k$-module of {\em loops} or cocharacters, i.e.,
$\ksep$-homomorphisms 
${\Gm}\ra T_{\sep}$.  
There is a natural pairing $T^\star\times T_\star\ra \Z$ that enables
us to identify $T_\star$ with the dual $\Hom(\Lambda,\Z)$ of $T^\star$, 
i.e., with the lattice of co-roots, 
denoted by $\coroot$ (see~\cite[VI.1.1, Prop.~2]{Bou:g4}). 
We fix a set of simple roots $\D = \{ \alpha_1, \alpha_2, \ldots,
\alpha_\ell \}$ of 
$\Phi$; the coroots $\ach_1, \ach_2, \ldots, \ach_\ell$ are a set of simple
roots of 
the dual root system $\Pch$ \cite[VI.1.5, Rem.~5]{Bou:g4}.  
Given an element $\ach := \sum c_i \ach_i$ in $\coroot$, the
corresponding loop
$\Gm\ra T_{\sep}$ is 
\begin{equation} \label{St.tor}
t \mapsto \prod h_i(t^{c_i}),
\end{equation}
where $h_i \!: \Gm \ra T_{\sep}$ is the loop corresponding to $\ach_i$.  The map 
$\prod h_i \!: \Gm^{\times \ell} \ra T_{\sep}$ is an isomorphism \cite[p.~44, Cor.~(a)]{St}.  
We have the identity:
\begin{equation} \label{St.eq}
h_\ach(t) x_\beta(u) h_\ach(t)^{-1} = x_\beta(t^{(\beta, \ach)} u) \qquad (\beta \in \Phi, \ach \in \Pch).
\end{equation} 
\end{borel}

\begin{borel}{Definition of $G'$} \label{Gpdef}
Assume that $G$ satisfies \eqref{cond.ind}.  Write $\coweight$ for the weight lattice of 
the dual root system; it has a basis, dual to the basis $\D$ of $\La_r$, and 
we write $\och_j$ for the ``co-weight" such that
$(\och_j, \alpha_i) = \delta_{ij}$ (Kronecker delta). Write $\Mbar$ for the sublattice of 
$\coweight$ generated by the $\och_j$ for $\alpha_j \in \D_r$, and $M$ for the 
intersection $\Mbar \cap \coroot$.  It follows from \cite[6.7, 6.9]{BoTi} that $\G_k$ 
acts trivially on $\Mbar$ --- hence also on $M$ --- because $\D_r$ is
pointwise fixed by the Galois action.  
(Compare \cite[5.2]{M:norm}.)  Consequently, the loops in $M$ are 
$k$-homomorphisms $\Gm \ra T$, and they generate a $k$-split torus $S'$ of rank $| \D_r |$ in $G$.

We define $G'$ to be the derived subgroup of the centralizer in $G$ of $S'$.

Over $\ksep$, we can describe $G'$ concretely.  It follows from \eqref{St.eq} that $G'$ is generated by the images of the $x_\alpha$'s, where $\alpha$ varies over the roots in $\Phi$ whose support does not meet $\D_r$.  In particular, $G'$ is semisimple and even simply connected by \cite[5.4b]{SpSt}.  The Dynkin diagram of $G'$ is obtained by deleting the vertices $\D_r$ from the Dynkin diagram of $G$.
The intersection of $G'$ with the maximal torus $T$ is, again over $\ksep$, the image of $\prod_{\alpha_i \not\in \D_r} h_i$.
\end{borel}

We are going to compute $t_{R,G}$
by reducing to this subgroup $G'$. One reason why this is possible is
the following: 

\begin{prop} \label{center} 
The center $Z$ of $G$ is contained in $G'$. 
\end{prop} 

\begin{proof}
It suffices to check this over an algebraic closure of $k$, so we may
assume that $G'$ and $G$ are split.  It follows from equation \eqref{St.eq} that the center of $G$ is the subgroup of $T$ cut out by the roots, i.e., it is the intersection
\[
Z(G) = \cap_{\alpha \in \Phi} \ker \alpha.
\]
Of course, we may replace the condition ``$\alpha \in \Phi$" with ``$\alpha$ is in the root lattice".  In particular, the fundamental weight $\omega_j$ corresponding to $\alpha_j \in \D_r$ belongs to the root lattice, so $Z(G)$ is contained in the kernel of $\omega_j$.  Recall that $\omega_j$ is given by the formula
\[
\omega_j\left(\prod\nolimits_i h_i(t_i) \right) = t_j.
\]
Therefore, $Z(G)$ is a subgroup of $\prod_{\alpha_i \not\in \Delta_r} h_i(\Gm)$.  But this is the maximal torus of $G'$.
\end{proof}

\begin{borel*} \label{red.lem}
As $G'$ is simply connected, it is of the form $G'_1 \times G'_2 \times \cdots \times G'_s$, where each $G'_i$ is simply connected and absolutely simple with Dynkin diagram a connected component of $\D \setminus \D_r$.

\begin{lem*}
An element $a \in Z^1(k, \prod_i G'_i)$ can be written as $\prod_i
a_i$ for $a_i \in Z^1(k, G'_i)$, and the Rost invariant $r_G(a)$ equals $\sum_i m_i r_{G'_i}(a_i)$, 
where $m_i$ is the Rost multiplier of the inclusion $G'_i \subset G$.
\end{lem*}

\begin{proof}
We prove this by induction on the number $r$ of $a_i$'s such that $a_i$
is not the zero cocycle in $H^1(k, G'_i)$.  
The case $r = 1$ is trivial.  
We treat the general case.  
Let $i$ be such that $a_i$ is not zero.  
Consider the group $G_{a_i}$ obtained by twisting $G$ by $a_i$; 
since $G'_i$ commutes with $G'_j$ for $j \ne i$, the twisted group
contains $\prod_{j \ne i} G'_j$ 
in an obvious way.  We find a diagram
\[
\begin{CD}
H^1(k, \prod_{j \ne i} G_j) @>>> H^1(k, G_{a_i}) @>{r_{G_{a_i}}}>> H^3(k, \QZt) \\
@. @V{\tau_{a_i}}V{\cong}V @VV{? + r_G(a_i)}V \\
@. H^1(k, G) @>{r_G}>>H^3(k, \QZt)
\end{CD}
\]
where $\tau_{a_i}$ denotes the twisting map.  
The diagram commutes as in Remark \ref{Deven.rmk}.
Starting with $\prod_{j \ne i} a_j$ in the upper left, 
we find $\sum_{j \ne i} m_j r_{G'_j}(a_j)$ in the upper right by the
induction 
hypothesis and $a$ in $H^1(k, G)$ in the lower left.  
The commutativity of the diagram gives the desired equation.
\end{proof}
\end{borel*}

\begin{borel*} \label{reduce.rost}
Continue the notation of \ref{red.lem}.  We are mainly interested in the case where $G$ is ``simply laced'',
i.e., where all roots have the same length.  
In that case, the inclusions $G'_i \subset G$ all have Rost multiplier
one---see e.g.\ \cite[2.2]{G:rinv} 
or \cite[7.9.2]{MG}---and the formula from Lemma \ref{red.lem} simply says:
\begin{equation} \label{red.simple}
r_G(a) = \sum_i r_{G'_i}(a_i).
\end{equation}
\end{borel*}

\section{The center and the torus} \label{center.1} 

For the duration of this section, the field $k$ is arbitrary, with no restriction on the characteristic.  We compute how the scheme-theoretic center $Z$ of a semisimple simply connected group $G$ sits inside a fixed maximal torus, in terms of the root system.  We do this using an exact sequence from \cite{M:norm}, which we first recall. 

\borel{Isogenies of tori}\label{vos}
Let $T \ra \Tbar$ be a $k$-isogeny of tori with kernel $Z$, 
i.e., $Z$ is finite and the sequence
\[
\begin{CD}
1 @>>> Z @>>> T @>\pi>> \Tbar @>>> 1
\end{CD}
\]
is exact.  
As before, $T_\star$ (resp. $\Tbar_\star$) is the $\Gamma_k$-module of {\em loops} or cocharacters, i.e.,
$\ksep$ homomorphisms 
${\Gm}\ra T_{\sep}$ (resp. $\Gm\ra \Tbar_{\sep}$).

\begin{prop} \label{vos.prop} 
There is a map $z$ such that
\[
\begin{CD}
1 @>>> T_\star @>{\pi_\star}>> \Tbar_\star @>z>> \varinjlim_n \Hom_{\ksep}(\mmu{n}, Z_{\sep}) @>>> 1
\end{CD}
\]
is an exact sequence of $\Gamma_k$-modules.
\end{prop}

This is sketched in \S1 of Merkurjev's paper (loc.cit.).  
We give an explicit proof, including a precise description of $z$
which will be used later.  

We first prove the following lemma. 
\begin{lem} \label{lift}
The image of every $k$-homomorphism from $\mmu{n}$ to a split torus
$S$ is contained in a 
rank 1 subtorus of $S$.
\end{lem}

\begin{proof}
We may assume that $\mmu{n}$ injects into $S$.  Dualizing, such an
injection corresponds 
to a surjection $\pi \!: \Z^r \ra \Zm{n}$, where $r$ is the rank of
$S$. By the fundamental 
theorem for finitely generated abelian groups, there is a basis $b_1,
\ldots, b_r$ of $\Z^r$ such 
that the kernel of $\pi$ has basis $b_1, b_2, \ldots, b_{r-1}, nb_r$.
Projection on the 
$b_r$-coordinate defines a map $\pi_r \!: \Z^r \ra \Z$, and $\pi$
factors through $\pi_r$.  
Dualizing again gives the lemma.
\end{proof}

\begin{proof}[Proof of Prop.~\ref{vos.prop}]
Since $\pi$ has finite kernel, the induced map $\pi_\star$ is injective.  
We use it to identify $T_\star$ with a submodule of $\Tbar_\star$.  Since the
tori have the same rank, 
the quotient $\Tbar_\star / T_\star$ is finite.

We now define $z$.  For $\omega \in \Tbar_\star$, let $n$ be the smallest
positive integer 
such that $n\omega$ belongs to $T_\star$.  We define $z_\omega \!: \mmu{n}
\ra T_{\sep}$ 
to be the restriction of the loop $n\omega \!: \Gm \ra T_{\sep}$ to
${\mmu{n}}$.  
Every character $\chi$ of $\Tbar$ maps $\Tbar_\star \ra \Z$, hence
$\chi(n\omega)$ is in $n\Z$ 
and the composition
\[
\begin{CD}
\mmu{n} @>{z_\omega}>> T @>>> \Tbar @>{\chi}>> \Gm
\end{CD}
\]
is 1.  Since subgroups of $T$ are cut out by characters, the image of
$z_\omega$ belongs to $Z$.  
It is an exercise to verify that the sequence displayed in the lemma
is exact at $\Tbar_\star$.  
Note that the map $z$ is $\Gamma_k$-equivariant.  

It remains to prove that $z$ is surjective.  Let $\phi \!: \mmu{n} \ra
Z_{\sep}$ be a homomorphism 
defined over $\ksep$.  By Lemma \ref{lift}, there is a rank 1 subtorus
$S_{\sep}$ of $T_{\sep}$ such that 
the image of $\phi$ is contained in $S_{\sep}$; let $\rho$ be a loop
generating $S_\star$.  
By hypothesis, $\chi(\pi \rho)$ is in $n\Z$ for every character $\chi
\in \Tbar^\star$, 
so there is some $\omega \in \Tbar_\star$ such that $\pi \rho = n\omega$.  
Then $\phi$ and $z_\omega$ agree in $\varinjlim \Hom(\mmu{n}, Z_{\sep})$.
\end{proof}

In our situation (with the notations from \ref{genandrel} and
\ref{Gpdef}), the center $Z$ of $G$ is contained in the maximal torus $T$ and the sequence 
\[
\begin{CD} 
1 @>>> Z @>>> T @>>> \Tbar @>>> 1 
\end{CD}
\]
is exact. 
For the same reason that $T_\star$ is identified with the coroot lattice $\coroot$, the lattice $\Tbar_\star$ is identified with $\coweight$. 
Hence the exact sequence given by Prop.~\ref{vos.prop} can be
re-written as 
\begin{equation}\label{coroot}
\begin{CD} 
1 @>>> \coroot @>>> \coweight @>>> \varinjlim_n \Hom_{\ksep}(\mmu{n},
Z_{\sep}) @>>> 1
\end{CD}
\end{equation} 

We close this section with an application of this exact sequence.  We will only use it for groups of type $B$ and $C$, but we include it because of independent interest.  Write $\D_c$ for the set of $\alpha_j \in \D$ such that $\och_j$ (as defined in \ref{Gpdef}) is a minuscule weight for the dual root system $\Pch$.  Recall \cite[\S{VI.1}, Exercise 24]{Bou:g4} that the minuscule weights are the minimal nonzero dominant weights.

\begin{cor} \label{vanish}
If $G$ is of inner type and $\D_c$ does not meet $\D_0$, then the center $Z$ of $G$ is contained in the maximal $k$-split subtorus $S$ of $T$.
\end{cor}

The condition on $\D_c$ is the same as saying: every element of $\D_c$ is circled in the index of $G$.


\renewcommand{\theenumi}{\roman{enumi}}
\begin{rmk} \label{vanish.rmk}
In the situation of the corollary, $H^1(k, S)$ is zero, so the natural map $H^1_{\mathrm{fppf}}(k, Z) \ra H^1(k, G)$ is zero.  (If $Z$ is smooth---as in the rest of this paper---the group $H^1_{\mathrm{fppf}}(k, Z)$ agrees with the Galois cohomology group $H^1(k, Z)$.)  This has two useful consequences:
\begin{enumerate}
\item The composition \eqref{hom} is zero.
\item The connecting homomorphism $(G/Z)(k) \ra H_{\mathrm{fppf}}^1(k, Z)$ is surjective.
\end{enumerate}
\end{rmk}

\begin{proof}[Proof of Cor.~\ref{vanish}]
Since the minuscule weights generate $\coweight /  \coroot$
\cite[\S{VI.2}, Exercise 5]{Bou:g4},  $Z$ is generated by the images
of $z_{\och_j}$ for $\alpha_j \in \D_c$.  Using the same argument as
in the proof of Prop.~\ref{vos.prop}, we see that $Z$ is contained in the torus $Q$ (defined over $\ks$) corresponding to the $n_j \och_j$ for $\alpha_j \in \D_c$, where $n_j$ is the smallest natural number such that $n_j \och_j$ is in $\coroot$.

Since $G$ has inner type, the maximal $k$-split torus $S$ in $T$ is the intersection of the kernels of $\alpha \in \D_0$ by \cite[6.7, 6.9]{BoTi}.  For such an $\alpha$ and for $\alpha_j \in \D_c$, we have $\alpha \ne \alpha_j$ by hypothesis, so the inner product $(\alpha, \och_j)$ is zero.  We conclude that the torus $Q$ is contained in $S$.
\end{proof}

\section{Type $B$} \label{B.section}

Before proving the main theorem for exceptional groups, 
we show that our method gives a new and very short proof for groups of
type $B_\ell$ ($\ell \ge 2$), i.e., that composition \eqref{hom} is zero.  
By \ref{Titsindex.prop}, it is enough to prove this for groups
having Tits index satisfying~\eqref{cond.ind}.  Consulting the tables in \cite{Bou:g4}, we find that 
\[
\D_r = \{ \alpha_1, \alpha_2, \ldots, \alpha_{\ell - 1} \}.  
\]
The dual root system is $C_\ell$, and its only minuscule weight is $\och_1$, i.e., 
\[
\D_c = \{ \alpha_1 \}.
\]
By \eqref{cond.ind}, the vertex $\alpha_1$ is circled in the index of $G$, so by Remark \ref{vanish.rmk}(i) the composition \eqref{hom} is zero.

Note that groups of type $B_\ell$ satisfying the hypotheses of Cor.~\ref{vanish} are those of the form $\Spin(q)$ where $q$ is a $(2 \ell +1)$-dimensional \emph{isotropic} quadratic form.  For these groups, Remark \ref{vanish.rmk}(ii) says: The spinor norm map $SO(q)(k) \ra \ksq$ is surjective.

\section{Concrete description of the center} \label{center.2}

We now give a concrete description of the center $Z$ inside a simply connected group $G$ in terms of the generators and relations for $G$ as in \ref{genandrel}; consequently we work over the separably closed field $\ks$.  (But see \ref{center.inner} below.)
Note that, by exactness of 
sequence \eqref{coroot}, it is enough to compute $z_{\och}$ for $\och$ 
belonging to a set of representatives in $\coweight$ of generators of the
quotient group $\coweight/\coroot$.

\begin{eg}[$E_7$] \label{E7.eg}
The center of a simply connected group $G$ of type $E_7$ is
isomorphic to $\mmu2$.  
As all roots have the same length, we can normalize the Weyl-invariant
inner product 
so that all roots have length 2; this identifies $\Phi$ with the
inverse system $\Pch$.  
The fundamental weight
\[
\omega_7 = \alpha_1 + \frac32 \alpha_2 + 2 \alpha_3 + 3\alpha_4 
+ \frac52 \alpha_5 + 2\alpha_6 + \frac32 \alpha_7
\]
is not in the root lattice and so maps to the nonidentity element in
$\La / \La_r$.  (We systematically number 
roots as in Bourbaki; for a diagram see \S\ref{E7.sec}.)  The
corresponding map $z_{\omega_7} \!: \mmu2 \ra Z_{\sep}$ is 
given by
\[
z_{\omega_7}(-1) = h_{2\omega_7}(-1) = h_2(-1) h_5(-1) h_7(-1).
\]

Here we can see Prop.~\ref{center} explicitly: $\D\setminus\D_r$ is $\{\alpha_2,\alpha_5,\alpha_7\}$,
and the image of this morphism is contained in $G'$.
\end{eg}

\begin{eg}[$E_6$] \label{E6.eg} 
Let now $G$ be of type $E_6$. The center of $G$ is isomorphic to
$\mmu{3}$. Again we can identify $\Phi$ and $\Pch$. 
The fundamental weight 
\[
\omega_1 = \frac43 \alpha_1 + \alpha_2 + \frac53 \alpha_3 + 2\alpha_4 + \frac43 \alpha_5 + \frac23 \alpha_6
\]
is not in in the root lattice; so its class generates
$\Lambda/\Lambda_r$. 
The corresponding isomorphism $z_{\omega_1} \!: \mmu3 \ra Z_{\sep}$ is defined by 
\begin{equation} \label{E6.equation}
z_{\omega_1}(\zeta) = h_{3\omega_1}(\zeta) = h_1(\zeta) \, h_3(\zeta^2) \, h_5(\zeta) \, h_6(\zeta^2).
\end{equation}
\end{eg} 


\begin{eg}[$C_\ell$] \label{Cn.eg} 
Let now $G$ be of type $C_\ell$; the center $Z$ is isomorphic to
$\mmu2$. In this case, 
we must be somewhat more careful because the root system has roots of
different lengths.  
The inverse system $\Pch$ is $B_\ell$.  We write $\och_\ell$ for the
fundamental weight corresponding to the root 
$\ach_\ell$ in the inverse system, i.e.,
\[
\och_\ell = \frac 12 (\ach_1 + 2\ach_2 + \cdots + \ell \ach_\ell).
\]
The corresponding isomorphism $\mmu{2}\ra Z_{\sep}$ is given by 
\begin{equation} \label{C.eq} 
z_{\och_\ell}(-1) = \begin{cases}
h_{2\och_\ell}(-1) = h_1(-1) h_3(-1) \cdots h_{\ell-1}(-1) & \text{if $\ell$ is even,} \\
h_{2\och_\ell}(-1) = h_1(-1) h_3(-1) \cdots h_{\ell}(-1) & \text{if $\ell$ is odd.}
\end{cases}
\end{equation} 
\end{eg} 

\begin{eg}[$D_\even$] \label{Dn.eg}
For $G$ of inner type $D_\ell$ with $\ell$ even, the center $Z$ is isomorphic to
$\mmu2 \times \mmu2$.  
As for $E_7$, all roots have the same length and we identify the root
system with the inverse system.  
The fundamental weights
\begin{gather*}
\omega_{\ell-1} = \frac12 \alpha_1 + \alpha_2 + \frac32 \alpha_3 + \cdots
+ 
\frac{\ell-2}{2}\alpha_{\ell-2} + \frac{\ell}4 \alpha_{\ell-1} + \frac{\ell-2}4 \alpha_\ell \\
\omega_{\ell} = \frac12 \alpha_1 + \alpha_2 + \frac32 \alpha_3 + \cdots + 
\frac{\ell-2}{2}\alpha_{\ell-2} + \frac{\ell-2}4 \alpha_{\ell-1} + \frac{\ell}4 \alpha_\ell 
\end{gather*}
have two different residues in $\La/ \La_r$.  As this is the Klein
four-group, they generate it.  
It follows that the images of the corresponding maps
$z_{\omega_{\ell-1}}, z_{\omega_\ell} \!: \mmu2 \ra Z_{\sep}$ generate $Z_{\sep}$.  
Applying \ref{vos}, we find that $Z_{\sep}$ is generated by the images of the homomorphisms $z_0, z_1$ defined by
\begin{align} 
z_0(-1) &:= h_1(-1) h_3(-1) \cdots h_{\ell-3}(-1) \, h_\ell(-1) \notag \\
z_1(-1) &:= h_1(-1) h_3(-1) \cdots h_{\ell-3}(-1)\, h_{\ell-1}(-1) \label{Dn.maps}
\end{align}
Note that we have abandoned tying the weights with their corresponding
homomorphism, 
because this depends on the parity of $\ell/2$.
\end{eg}

\begin{eg}[$A_2$]\label{A2.eg} 
Suppose now that $G$ has type $A_2$, i.e., $G$ is isomorphic to $SL_3$ and its center is isomorphic to $\mmu3$.
Since the fundamental weight $\omega_1=\frac 23 \alpha_1+\frac 13
\alpha_2$ does not belong to the root lattice, we get an isomorphism 
$z_{\omega_1}:\,\mmu{3}\ra Z_{\sep}$ defined by 
\begin{equation}\label{A2.equation} 
z_{\omega_1}(\zeta)=h_1(\zeta^2)h_2(\zeta).
\end{equation}  
\end{eg}

\begin{borel*} \label{center.inner}
Let $G$ be a simply connected semisimple group of inner type over $k$.  Each element $\och$ of $\coweight$ defines an isomorphism
\[
z_\och \!: \mmu{n} \ra \mmu{n} \subset G
\]
over $\ksep$.  But the only morphisms $\mmu{n} \ra \mmu{n}$ are raising to a power, so $z_\och$ is $\G_k$-equivariant, and hence is defined over $k$ by \cite[AG.14.3]{Borel}.
\end{borel*}

\section{Type $E_7$} \label{E7.sec}

This section proves the main theorem \ref{mainthm.thm} for a group $G$ of type $E_7$.  
In this case, the center $Z$ is identified with $\mmu2$.  
We prove that $t_{R,G}$ coincides with the Tits class $t_G$, i.e. 
the composition \eqref{hom} is the cup product with $t_G$, 
where the cup product is induced by the obvious bilinear map 
$Z(\ks) \times Z(\ks) = \mmu2(\ks) \times \mmu2(\ks) \ra
\mmu2(\ks)^{\otimes 2}$. 
(The group $H^2(k, Z)$ is 2-torsion, so $t_{R,G}$ and $t_G$ generate the same subgroup of $H^2(k, Z)$ if and only if they are equal.)  

By \ref{Titsindex.prop}, it is enough to prove this for groups of type $E_7$ 
with Tits index satisfying \eqref{cond.ind}. 
Let $G$ be such a group, i.e. assume the index of $G$ is 
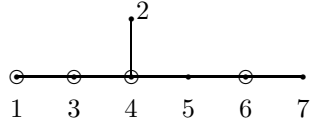
\setlength{\unitlength}{0.3in}
\[
\begin{picture}(7,3)
    \multiput(1,1)(1,0){6}{\circle*{\darkrad}}
    \put(3,2){\circle*{\darkrad}}

    \put(1,1){\line(1,0){5}}
    \put(3,2){\line(0,-1){1}}
    
    \put(1,0.3){\makebox(0,0.4)[b]{$1$}}
    \put(2,0.3){\makebox(0,0.4)[b]{$3$}}
    \put(3,0.3){\makebox(0,0.4)[b]{$4$}}
    \put(4,0.3){\makebox(0,0.4)[b]{$5$}}
    \put(5,0.3){\makebox(0,0.4)[b]{$6$}}
    \put(6,0.3){\makebox(0,0.4)[b]{$7$}}
    \put(3.2,2){\makebox(0,0.4)[b]{$2$}}
    
    \multiput(1,1)(1,0){3}{\circle{\lrad}}
    \put(5,1){\circle{\lrad}}

\end{picture}
\]
or possibly has more circled vertices.  (If the index has more
circles, 
then $G$ is split, see \cite[p.~57]{Ti:Cl}.)

Write $G'_i$ with $i = 2, 5, 7$ for the components of $G'$, where
$G'_i$ corresponds to the vertex 
$\alpha_i$ in the Dynkin diagram.  
Each is isomorphic to $SL_1(Q_i)$ for some quaternion algebra $Q_i$.  
Since the weights $\omega_2, \omega_5, \omega_7$ all have the same
image in $\La /\La_r$, 
the $Q_i$ are all isomorphic by \cite[p.~211]{Ti:R} and we simply write $Q$.  By \cite[loc. cit.]{Ti:R},
$[Q] \in H^2(k, \mmu2)$ is also the Tits class of $G$.
The center $Z$ is contained in the product $Z_2 \times Z_5 \times
Z_7$ of the centers of $G'_2, G'_5, G'_7$ 
by Prop.~\ref{center}.  
By the description of the center given in example \ref{E7.eg}, the map 
\[
\mmu2 = Z \ra Z_2 \times Z_5 \times Z_7 = \mmu2 \times \mmu2 \times
\mmu2
\] is 
$-1 \mapsto (-1, -1, -1)$.  (Even though Example \ref{E7.eg} treats $G$ over the separable closure, it still applies to our nonsplit $G$ over $k$ by \ref{center.inner}.)
Applying equation \ref{red.simple} and the type $A_1$ case of the main
theorem (see Remark~\ref{typeA.rmk}), 
we get that  the class $t_{R,G}$ is 
$3[Q] = [Q]$, which is the Tits class of $G$.

This concludes the proof of the main theorem for groups of type $E_7$.

\section{Type $C$} \label{C.section}

The argument of the previous section can be adapted to recover the
theorem for groups of type $C_\ell$.   (Alternatively, in the case where $\ell$ even, the same argument as for type $B$ in \S\ref{B.section} shows that composition \eqref{hom} is zero.)
Since the result is not new for type $C$, we only briefly sketch the proof. 
We want to show that $t_{R,G}$ is $0$ if $\ell$ is even and $t_G$ if $\ell$
is odd. 
Recall that a group of type $C_\ell$ has Dynkin diagram
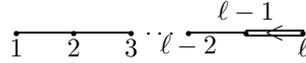
\[
\begin{picture}(7,2)
    \put(1,1){\line(1,0){2}} 
    \put(4,1){\line(1,0){1}}
    \put(5,1.04){\line(1,0){1}}
    \put(5,0.96){\line(1,0){1}}
 
    \put(1,0.6){\makebox(0,0.4)[b]{$1$}}
    \put(2,0.6){\makebox(0,0.4)[b]{$2$}}
    \put(3,0.6){\makebox(0,0.4)[b]{$3$}}
    \put(4,0.6){\makebox(0,0.4)[b]{$\ell-2$}}
    \put(5,1.2){\makebox(0,0.4)[b]{$\ell-1$}}
    \put(6, 0.6){\makebox(0,0.4)[b]{$\ell$}}
    
    \put(5.5,0.87){\makebox(0,0.4)[b]{$<$}}

    \multiput(3.3,1)(0.2,0){3}{\circle*{0.01}}


    \multiput(1,1)(1,0){6}{\circle*{\darkrad}}

\end{picture}
\]
As before, it is enough to prove it for groups whose Tits
index satisfy condition~\eqref{cond.ind}, i.e., such that $\alpha_i$ is circled for $i$ even.

As before, we write $G'_i$ for the component
of $G'$ corresponding to the vertex $\alpha_i$, for $1\leq i\leq n$, $i$
odd. We may identify all those with $SL(Q)$ for some quaternion
algebra $Q$. 
The description of the center given in
Example~\ref{Cn.eg} shows that $Z$ maps to the product of the centers
of the $G'_i$ by  $(-1)\mapsto(-1,\dots,-1)$. 

It remains to apply Lemma~\ref{red.lem}, which asserts 
$r_G(a)=\sum_i m_i r_{G'_i}(a_i)$.
Note that we do have to care about Rost multipliers here, since groups
of type $C$ are not simply laced. 
If $1\leq i\leq n-1$, $i$ odd, then the root $\alpha_i$ is short in
$C_\ell$. So the corresponding co-root $\ach_i$ is long and the inclusion 
$G'_i\subset G$ has Rost multiplier $2$, see e.g.\ \cite[2.2]{G:rinv} 
or \cite[7.9.2]{MG}.  Since the Rost invariant 
has values in $H^3(k,\mmut{2})$ which is $2$-torsion, 
we get $r_G(a)=0$ if $\ell$ is even, and 
$r_G(a)=r_{G'_\ell}(a_\ell)$ if $\ell$ is odd. 
In this last case, applying the formula for groups of type $A$, we get 
that $t_{R,G}=[Q]$, which is the Tits class of $G$. 

\section{Type $E_6$} \label{E6.sec}

In this section, we prove the theorem for a group $G$ of type
$E_6$.

\borel{Groups of inner type $E_6$} 
Let us first assume $G$ is of inner type, that is $\Gamma_k$ acts
trivially on the Dynkin diagram. 
The center $Z$ is identified with $\mmu{3}$ and we want to
prove that $t_{R,G}$ and $t_G$ generate the same subgroup of $H^2(k, Z)$.
Again by \ref{Titsindex.prop}, it suffices to consider
groups whose Tits index satisfies condition~\eqref{cond.ind}. Let $G$ be
such a group, i.e. we assume the Tits index of $G$ is 
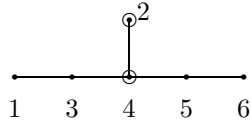
\[
\begin{picture}(7,3)
    \multiput(1,1)(1,0){5}{\circle*{\darkrad}}
    \put(3,2){\circle*{\darkrad}}

    \put(1,1){\line(1,0){4}}
    \put(3,2){\line(0,-1){1}}
    
    \put(1,0.3){\makebox(0,0.4)[b]{$1$}}
    \put(2,0.3){\makebox(0,0.4)[b]{$3$}}
    \put(3,0.3){\makebox(0,0.4)[b]{$4$}}
    \put(4,0.3){\makebox(0,0.4)[b]{$5$}}
    \put(5,0.3){\makebox(0,0.4)[b]{$6$}}
    \put(3.25,2){\makebox(0,0.4)[b]{$2$}}
    
    \put(3,1){\circle{\lrad}}
    \put(3,2){\circle{\lrad}}
\end{picture}
\]
or has more circled vertices.  
(As in the type $E_7$ case in \S\ref{E7.sec}, if more vertices are
circled, then $G$ is split.)

Write $G'_1$ and $G'_5$ for the components of $G'$ corresponding to
the subdiagrams 
with vertices 1, 3 and 5, 6 respectively.  
They are of type $A_2$. 
Moreover, since the fundamental weights $\omega_1$ and $\omega_5$ have
the same class in $\Lambda/\Lambda_r$, you may identify $G'_1$ and
$G'_5$ with $SL(D)$ for some central
simple algebra $D$ of degree $3$ over $k$.  The center of $SL(D)$ is naturally identified with $\mmu3$, and the composition
\[
H^1(k, \mmu3) \ra H^1(k, SL(D)) = H^1(k, G'_i) \xrightarrow{r_{G'_i}} H^3(k, \QZt)
\]
is the cup product with $m[D]$ for $m = 1$ or 2 by Example \ref{typeA.rmk}.

Comparing the description of the center of $G$ and the $G'_i$ given by 
formulas~\ref{E6.equation} and~\ref{A2.equation}, one sees that the
map $Z\ra G'_1\times G'_5$ is $\zeta \mapsto (\zeta^2,\zeta^2)$.  
Hence the induced map $H^1(k,Z)\ra H^1(k,Z'_1)\times H^2(k,Z'_5)$ is
given by $a\mapsto (2a,2a)$. 
Applying Equation \eqref{red.simple}, we see that the composition \eqref{hom} is 
\[
a \mapsto 4a \cp m[D] \quad \in H^3(k, \mud3).
\]
But $4a$ equals $a$ because the exponent of $Z$ is $3$.  Therefore, $t_{R,G}$ equals $m[D]$, i.e., generates the same subgroup of $H^2(k, Z)$ as $[D]$, which is $t_G$ \cite[p.~211]{Ti:R}.

\borel{Groups of outer $E_6$ type}
Suppose now that the group $G$ is of outer type, 
i.e., the absolute Galois group acts nontrivially on the Dynkin
diagram of $G$.  
The kernel of the homomorphism $\Gamma_k \ra \aut(\D)=\Z/2$ 
defines a quadratic extension $L$ of $k$, over which the group becomes of inner type. 
Its center $Z$ is isomorphic to ${\mmu{3}}_{[L]}$. 

By the inner type case, we know that $t_{R,G} - m\,t_G$ is killed by $L$. 
Since it lives in
the group $H^2(k,Z)$ which is of exponent $3$, while $L$
is quadratic over $k$, a 
restriction/corestriction argument shows that $t_{R,G} - m\,t_G$ is zero in $H^2(k, Z)$.
This concludes the proof for groups of type $E_6$. 

\section{Type $D_\ell$ for $\ell$ even} \label{D.section} 

We now prove 
the main theorem for groups of type $D_\ell$ for $\ell$ even,
including trialitarian $D_4$. 

\begin{borel}{Groups of inner type $D_\ell$} 
As before, we first assume that $G$ is of inner type, so that $Z$ is
identified with $\mmu{2}\times \mmu{2}$. 
We want to prove that the composition~(\ref{hom}) is the cup product
with the Tits class, where the cup product is induced by the bilinear
map given in~(\ref{Deven.cup}). 
By \ref{Titsindex.prop}, it is enough to consider groups with Tits
index satisfying condition~\eqref{cond.ind}, i.e., with Dynkin diagram
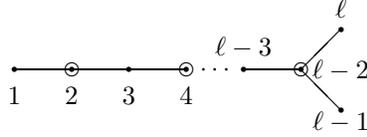
\[
\begin{picture}(9,2)
    \multiput(1,1)(1,0){3}{\line(1,0){1}}
    \put(5,1){\line(1,0){1}}
    
    \put(1,0.4){\makebox(0,0.4)[b]{$1$}}  
    \put(2,0.4){\makebox(0,0.4)[b]{$2$}}
    \put(3,0.4){\makebox(0,0.4)[b]{$3$}}
    \put(4,0.4){\makebox(0,0.4)[b]{$4$}}
    \put(5,1.2){\makebox(0,0.4)[b]{$\ell-3$}}
    \put(6.2,0.8){\makebox(0,0.4)[l]{$\ell-2$}}

    \multiput(4.3,1)(0.2,0){3}{\circle*{0.01}}

    \put(6,1){\line(1,-1){0.7}}
    \put(6,1){\line(1,1){0.7}}

    \multiput(1,1)(1,0){6}{\circle*{\darkrad}}
    
    \put(2,1){\circle{\lrad}}
    \put(4,1){\circle{\lrad}}
    \put(6,1){\circle{\lrad}}

    \put(6.7,0.3){\circle*{\darkrad}}
    \put(6.7,1.7){\circle*{\darkrad}}
    
    \put(6.7,-0.1){\makebox(0,0.4)[b]{$\ell-1$}}
    \put(6.7,1.9){\makebox(0,0.4)[b]{$\ell$}}
    
\end{picture}
\]
or having more circled vertices. 

We write $G'_i$, for $i=1,3,\dots,\ell-3,\ell-1$ and $\ell$ for the
component of $G'$ corresponding to the vertex $\alpha_i$. 
They are of type $A_1$, and may be identified with 
$SL(Q_i)$, for some quaternion algebras $Q_i$. 
Each has center $Z'_i$ isomorphic to $\mmu{2}$. 
Moreover, since $\omega_1,\dots,\omega_{\ell-3}$ and
$\omega_{\ell-1}+\omega_\ell$ have the same class in
$\Lambda/\Lambda_r$, we may assume $Q_1=\dots =Q_{\ell-3}=Q$, 
for some quaternion algebra $Q$ satisfying $[Q]=[Q_{\ell-1}]+[Q_\ell]$. 

Since $G$ is of inner type, the homomorphisms $z_0, z_1$ from Example
\ref{Dn.eg} are 
$k$-defined even though $G$ need not be split.  
Hence, we get a map
\[
(z_0,z_1) \!:
\mu_2\times\mu_2\ra Z \subset Z'_1\times\dots\times Z'_{\ell-3}\times
Z'_{\ell-1}\times Z'_\ell.
\] 
The induced map $H^1(k,Z)\ra \prod_iH^1(k,Z'_i)$ is given by 
\[
(a_0,a_1)\mapsto (a_0+a_1,\dots,a_0+a_1,a_1,a_0).
\]
Applying Equation \eqref{red.simple}, and the formula for groups of
type $A$ (see~\ref{typeA.rmk}), we get 
$r_G(a_0, a_1)=\frac{\ell-2}2(a_0+a_1)\cp [Q]+ a_1\cp
[Q_{\ell-1}]+a_0\cp[Q_\ell],$
that is 
\begin{equation} \label{Deven.1}
r_G(a_0,a_1)=\begin{cases}
a_0\cp [Q_{\ell-1}]+a_1\cp [Q_{\ell}] &\text{if $\ell \equiv 0 \pmod{4}$}\\
a_0\cp [Q_\ell]+a_1\cp [Q_{\ell-1}] &\text{if $\ell\equiv 2 \pmod{4}$}
\end{cases}
\end{equation}

On the other hand, we have:
\begin{equation}\label{Deven.2}
\omega_{\ell - i}(a_0, a_1) = a_i \quad \text{and} \quad \omega_{\ell - i}(t_G) = [Q_{\ell - i}]
\end{equation}
for $i = 0, 1$.
\end{borel}
Combining \eqref{Deven.1} and \eqref{Deven.2} with the definition of the cup product in \eqref{Deven.cup}, we conclude that composition \eqref{hom} equals the cup product with the Tits class $t_G$.

\begin{borel}{Groups of outer $D_\ell$ type}
For the groups of type $^2\!D_\ell$, one can reduce to the
$^1\!D_\ell$ case as in \cite[p.~817]{MPT}.   Alternatively, one can apply the same method as for groups of inner type.  In this case, $G'$ is a direct product of $(\ell - 2)/2$ copies of $SL(Q)$ and a transfer $R_{L/k} SL(Q')$, where $Q'$ is a quaternion algebra over a quadratic extension $L$ of $k$.  To apply \eqref{red.simple}, one needs to know the Rost invariant of $R_{L/k}SL(Q')$, which is specified in \cite[9.8]{MG}.

It remains only to treat the trialitarian groups, i.e., the groups $G$
of type $\iiiD$ or $\viD$.  Recall what this means.  The automorphism group of the Dynkin diagram $\D$ of $G$ is the symmetric group on three letters, and the superscript 3 or 6 denotes the size of the image of homomorphism $\G_k \ra \aut(\D)$.  The kernel of this homomorphism fixes at least one separable cubic field extension of $k$; we pick one and call it $L$.  The group $G$ is of type
$\oD$ or $\iiD$ over $L$.  As the exponent of $Z$ (i.e., 2) is
relatively prime to the dimension of $L/k$ (i.e., 3), 
the proof may be completed using some restriction/corestriction
argument as for groups of outer type $E_6$.
\end{borel}

\section{Application to groups of type $E_7$} \label{E7.app}

We now apply the main theorem to prove a result about groups of type $E_7$.  A simply connected group of type $E_7$ can be described in terms of ``gifts" as in \cite{G:e7}.  A gift is a triple $(A, \s, \pi)$ where $A$ is a central simple $k$-algebra of degree 56, $\s$ is a symplectic involution on $A$, and $\pi \!: A \ra A$ is $k$-linear and satisfies certain axioms, see \cite[Def.~3.2]{G:e7} for a precise statement.  The Brauer class of $A$ is the Tits class $t_G$.  The $k$-points of $G$ are the elements $a \in A^\times$ such that $\s(a) a = 1$ and $\Int(a)\, \pi = \pi \Int(a)$, where $\Int(a)$ denotes the automorphism $x \mapsto axa^{-1}$ of $A$.

The associated adjoint group $\Gbar$ is a subgroup of $PGL(A)$; it has $k$-points $\Int(a)$ for $a \in A^\times$ such that $\s(a) a$ is in $\kx$ and $\Int(a) \,\pi = \pi \Int(a)$.  The exact sequence
\[
\begin{CD}
1 @>>> \mmu2 @>>> G @>>> \Gbar @>>> 1
\end{CD}
\]
gives a connecting homomorphism
\[
\delta \!: \Gbar(k) \ra H^1(k, \mmu2) = \ksq.
\]
As in \cite[\S31]{KMRT}, one finds
\[
\delta(\Int(a)) = \s(a) a \in \ksq.
\]
If $\la \kxsq$ is in the image of $\delta$, then its image in $H^1(k, G)$ is zero and the composition \eqref{hom} sends $\la$ to zero.  The main theorem gives:
\begin{cor}
If $\la \kxsq$ is in the image of $\delta$, then $(\la) \cp [A] = 0$.
\end{cor}
This result was previously observed in the case $k = \R$ in \cite[6.2]{G:e7}.

{\small \subsection*{Acknowledgements} The first author thanks Universit\'e Paris 13 for its hospitality while some of the work on this paper was performed.}


\providecommand{\bysame}{\leavevmode\hbox to3em{\hrulefill}\thinspace}
\providecommand{\MR}{\relax\ifhmode\unskip\space\fi MR }
\providecommand{\MRhref}[2]{%
  \href{http://www.ams.org/mathscinet-getitem?mr=#1}{#2}
}
\providecommand{\href}[2]{#2}

\end{document}